\documentclass[10pt,a4paper]{amsart}

\usepackage{amsmath,amssymb,amsfonts,amsbsy,amscd,epsfig,graphicx}
\usepackage{amsthm}
\usepackage{epstopdf}

\newcommand{\RR}{\Bbb R}

\newtheorem{theorem}{{Theorem}}[section]
\newtheorem{lemma}[theorem]{{Lemma}}
\newtheorem{fact}[theorem]{{Fact}}

\newtheorem{proposition}[theorem]{{Proposition}}
\newtheorem{remark}[theorem]{{Remark}}

\newtheorem{proposition.fondamentale} [theorem]{{Proposition
fondamentale}}

\newtheorem{definition}[theorem]{{Definition}}

\newtheorem{remarque.terminologique}[theorem]{{Remarque
terminologique}}
\newtheorem{quest.algeb}[theorem]{{Question alg\'ebrique}}
\newtheorem{conj.geom}[theorem]{{Conjecture g\'eom\'etrique}}

\newcommand{\epsw}[2] 
        {\begin{array}{c} \hspace{-1.3mm}
       \raisebox{-4pt}{\epsfig{figure=#1.eps,width=#2}}
       \hspace{-1.9mm}\end{array}}
\newcommand{\epsh}[2]  
        {\begin{array}{c} \hspace{-1.3mm}
       \raisebox{-4pt}{\epsfig{figure=#1.eps,height=#2}}
       \hspace{-1.9mm}\end{array}}

\begin{document}

 \title{Real and discrete holomorphy: Introduction to an algebraic approach}

\author{Sylvain Barr\'e} \email{Sylvain.Barre@univ-ubs.fr , www.univ-ubs.fr/lmam/barre} \address{Universit\'e de Bretagne-Sud, 
Centre de Recherche, B\^atiment  Yves Coppens 
Campus de Tohannic BP 573, 56 017 VANNES, France}

 \author{Abdelghani Zeghib} \email{Zeghib@umpa.ens-lyon.fr,  www.umpa.ens-lyon.fr/\^~zeghib} \address{CNRS, UMPA, ENS Lyon
46  all\'ee d'Italie, 
 69364 Lyon cedex 07,  FRANCE }

 \date{\today}
\maketitle

\begin{abstract}

We consider spaces for which
there is a notion of harmonicity for complex
valued functions defined on them. For instance,  this is the case of Riemannian
manifolds on one hand, and
(metric) graphs on the other hand. We observe that it is then possible to
define an ``amusing''
 notion of holomorphic functions on them, and show how rigid  it is in some
cases.
 \\ \\
{\sc Resum\'e.} 
On peut parler 
   d'harmonicit\'e des fonctions \`a valeurs complexes,  d\'efinies sur des vari\'et\'es riemanniennes, ou sur des graphes m\'etriques, ainsi que sur d'autres espaces plus g\'en\'eraux. 
Nous observons ici, que dans tous ces cas, on peut aussi parler d'holomorphie de ces fonctions. Nous d\'ecortiquons  le cas classique des vari\'et\'es riemanniennes, et la cas discret de certains graphes. Nous  montrons, en particulier,  une rigidit\'e dans le cas discret, qui le distingue du cas classique.

\end{abstract}

\section{Introduction}

\subsection{Harmonic and holomorphic}

 It is a usual practice to compare harmonicity with holomorphy
in any context where they are simultaneously defined.  In the classical
case of real and complex
valued mappings defined on open subsets of
${\Bbb C}$, the dictionary is realized, locally, by means of the notion of
conjugate part.

In the general case of hermitian complex manifolds, since one notion is
metric, and the other is
differentiable, the comparison is available, only if some compatibility
between the metric and
 the complex structure
is fulfilled. Roughly speaking this corresponds to the Kaehler property.

\subsubsection{General framework, Support of harmonicity.}

Classically, harmonic functions are defined by means of a metric on a
Riemannian manifold.
However, there is  larger class
of spaces where
one can talk about them.  A substantial class is that of (measurable) spaces
endowed
with  Markov processes.
As particular interesting cases,  we have  random walks on groups,
and
graphs endowed
with their natural ``simplicial'' Markov chains. Recall here that for a
function
$f$, its Laplacian $\Delta f$  is the function:
$$ \Delta f (x)=  \frac{1}{\nu(x)} \sum_{y \sim x} f(y) - f(x)$$
where $y \sim x$ means $y$ is adjacent to $x$, and $\nu(x)$ is the valency of
$x$.

$\bullet$ In this paper,  we will always deal with complex valued functions
$\phi: X \to {\Bbb C}$, where
$X$ is  {\it a Riemannian manifold} or {\it a graph}.


\subsubsection{Holomorphic functions}


 Even in the most classical  case, there are some and
in fact deep differences between harmonicity and holomorphy. Maybe, the
silliest one is that
the (usual) product
 of two holomorphic functions is holomorphic, but this is not the case of
harmonic functions!


This suggests to define holomorphic functions just by forcing invariance
under multiplication:

\begin{definition} \label{definition.holomorphic}

 Let  $\phi: X \to {\Bbb C}$, where
$X$ is  a Riemannian manifold or  a graph. We say that $\phi$ is holomorphic
if both
$\phi$ and its square $\phi^2$ are harmonic.

\end{definition}

Recall the classical formula for functions on Riemannian manifolds:
$$\Delta (f g) = f \Delta g + g \Delta f + \langle \nabla f, \nabla g
\rangle$$ (it is remarkable that this is
also valid in the case of graphs, for a natural definition of the gradient
$\nabla$).

Therefore, if $X = {\Bbb R}^2 $  ( $= {\Bbb C}$), $\phi = f+ \sqrt{-1}g$ is
holomorphic in
the sense of the definition above, iff, $\parallel \nabla f \parallel =
\parallel \nabla g \parallel$, and
$ \langle \nabla f, \nabla g \rangle = 0$. This exactly means that $\phi$
is conformal, or to use
a more precise terminology, it is {\it semi-conformal} which  means 
that $\phi$ may
have singular points (i.e.
where $d \phi =0$).

\subsubsection{Some remarks}
${}$

$\bullet$ A  mapping $\phi: {\Bbb C} \to {\Bbb C}$ is
holomorphic, according with our definition, iff $\phi$ is holomorphic or
anti-holomorphic, in the usual sense.

Actually, an individual anti-holomorphic mapping (on any complex manifold)
can be restored to become
holomorphic.  However, when considering all holomorphic and anti-holomorphic
mappings together, it is a nuisance
to separate them. Indeed that is exactly the counterpart of non-linearity
in the equation
$\Delta \phi = \Delta \phi^2 =0$.  In the general case unlike that of
complex manifolds, this does not split into two linear equations
corresponding to
 holomorphy and anti-holomorphy.

However, for the sake of simplicity, we will keep our term ``holomorphic''
as it is defined, since in fact, we will
always deal with individual functions.

$\bullet$
In the classical case of a mapping $\phi: {\Bbb C} \to {\Bbb C}$,
conformality of $\phi$, implies conformality and
hence harmonicity  for all powers $\phi^N$, for all integer $N$, that is,
for all $N$,
the equation $$ \Delta \phi= \ldots =  \Delta \phi^N = 0 \quad \quad (*)_N$$
is satisfied (this has an obvious  meaning for $N = \infty$). More generally:

\begin{theorem} \label{quadratic}
${}$

 $\bullet$ For any Riemannian manifold $X$,  we have,  $$(*)_2
\Longrightarrow (*)_\infty$$

That is, a holomorphic function satisfies all $(*)_N$ for any $N$.

$\bullet$ For any $N >2$, there is a graph $X$ and a holomorphic function
on it, which  satisfies $(*)_N$ but
which does not  satisfy $(*)_{N+1}$.

$\bullet$ For a graph $X$, of finite valency, there is a finite $m = m(X)$, 
such that any
function satisfying
$(*)_m$ is constant (and hence satisfies trivially  $(*)_\infty$).

\end{theorem}


\subsubsection{A f urther motivation : Dirichlet problem.}


 Consider the simplest Dirichlet problem for
  {\it bounded} harmonic functions   defined on the unit disc $ {\Bbb D} $
of ${\Bbb C}$. Classical harmonic analysis theory allows one to define a
boundary value isomorphism : $$\phi \in h^\infty({\Bbb D}) \to
\partial_\infty \phi \in L^\infty
 ( \partial{\Bbb D})$$
from  the Hardy space of bounded harmonic functions on the disc,  to the
space of bounded measurable
 functions on the circle.

The Hardy space $H^\infty({\Bbb D})$ of holomorphic (in the usual sense)
bounded functions, is a part
of $h^\infty({\Bbb D})$, but its image in $L^\infty(S^1)$ by the boundary
value mapping, is very hard to
explicit. In analytic words, this certainly reduces to an invariance by
Hilbert transform, but this
cannot help in  understanding how holomorphic functions are more
``regular'' than harmonic ones.
It is singularly suggesting to bring out a criterion, ensuring that a
function $\psi  \in
L^\infty(S^1)$ is the boundary value of a holomorphic function!

Therefore, our motivation from the Dirichlet problem viewpoint, is that we
are introducing here
an ``abstract'' notion of holomorphy, which may help to understand the true
nature of
holomorphic functions. On the other hand, this notion of holomorphy for
functions on
$X$, allows one to specify a class of bounded measurable functions on its
Poisson boundary, which are more regular than the others.


\subsubsection{A more general framework.}

Actually, one can take any ring $F$, and consider $F$-valued functions
defined on $X$.
In the case of graphs, it is meaningful to speak of harmonicity.  As for
Riemannian manifolds,
one just  needs  the  ring to be topological, since the Laplacian can be
defined via its infinitesimal mean value
property.

\begin{definition} Let $X$ be a Riemannian manifold or a graph, and $F$ a
topological ring.
A mapping $\phi: X \to F$ is called $F$-holomorphic, iff $\phi$ and
$\phi^2$ are
harmonic.
\end{definition}

We have just rewritten the same definition given in the case 
$F = {\Bbb C}$. Our
goal is
to point out some examples  that seem exciting, and deserve attention. We
mention
here the case when $F$ is, the quaternionic field (over the reals), and the
matrix algebras
$M_n({\Bbb R})$ or $M_n({\Bbb C})$...

\begin{remark} { \em
As an example, in another direction, consider the D'Alembertian operator on
${\Bbb R}^2$,
 $\Box \phi= \frac{\partial}{\partial x} \frac{ \partial }{ \partial y}
\phi  $, where
$\phi: {\Bbb R}^2 \to {\Bbb R}^2$.  See the second ${\Bbb R}^2$ as
 a ring (in fact an ${\Bbb R}$-algebra) which is the product of ${\Bbb R}$
by itself.  Suppose $\phi$ is a solution of $\Box \phi = \Box \phi^2 = 0$.
It then follows,
essentially (up to a switch of factors...), that $\phi (x, y) = (f(x),
g(y))$. This means
that $\phi$ is a conformal mapping of the Minkowski space
$({\Bbb R}^2, dxdy)$.
}
\end{remark}

\begin{remark} {\em The previous example suggests to define a wave equation
on graphs.
Recall at this stage, that  more generally than the canonical Markov chain
on a graph, there are Markov chains with arbitrary non-equiprobable
transition functions, or
equivalently, the graph has weighted edges. This is also equivalent to endow
the graph with an adapted metric, which is just given by these weights.

 Now, the idea is to allow  weights to take negative values.
Analogously to the term ``pseudo-Riemannian'', one may call such graphs
``pseudo-metric''
As a diffusion operator, one obtains
something like a D'Alembertian, with more complicated properties.  They
seem well adapted and connected
to the dynamics on the graph.  It is particularly interesting to define,
amongst pseudo-metric graphs,
a notion of Lorentzian ones: there is only one negative direction which
corresponds to the ``time''.

}

\end{remark}

\subsection{Harmonic morphisms.}
The discussion above corresponds essentially to our chronological
motivation in  bringing out this notion of
holomorphy.
We then looked for possible correlated notions in the geometric literature.
It was difficult to conclude if we weren't,  hopefully, intuitively guided
to  harmonic morphisms.   After,
having a look on this wide theme, we discovered that our holomorphic
functions on Riemannian
manifolds, are nothing  but harmonic morphisms with values in ${\Bbb R}^2$.
This is however
obscure, and not clear to extract, even, in standard celebrated references
in this theory. In any case, motivation
of harmonic morphisms are never stated as ours here.
 On the other hand, in the discrete case, that is for graphs, a theory of
harmonic morphisms
is not yet completely established.  Maybe, the problem comes from the choice,
to see Riemannian manifolds and graphs as objects of the same category or
not. Anyway, let's anticipate
and say that our holomorphic functions can not be seen as
harmonic morphisms from graphs to ${\Bbb R}^2$.

 \subsubsection{Riemannian case.} Recall that for Riemannian manifolds,  
a harmonic morphism is a mapping
which preserves sheaves
of local harmonic real functions. In other words, $\phi: X \to Y$ is a
harmonic morphism,
if for any real harmonic function $f$ defined on an open subset of $Y$,
$ f \circ \phi$ is harmonic.

 \subsubsection{Graph case.} Obviously in this definition, one can let $X$ and $Y$, to be, both, or just
one among then, a graph
(instead of a Riemannian manifold). In particular, one can speak of
harmonic morphisms
  $X \to {\Bbb C}$, where $X$ is a graph.

 \subsubsection{A difference.} Harmonic morphisms in the discrete case 
do not behave as nicely as
 in the classical case. 
 Here is one difference between the two situations, 
with respect to our notions of holomorphy.
\begin{theorem} \label{non.quadratic}
Let $\phi: X \to {\Bbb C}$ be  a mapping.
${}$

$\bullet$ If $X$ is a Riemannian manifold, then $\phi$ is holomorphic, iff,
$\phi$ is a harmonic morphism.

$\bullet$
There is a  graph $X$
having a  holomorphic function $\phi: X \to {\Bbb C}$  which is not a
harmonic morphism.

\end{theorem}

In fact, a non-constant  harmonic morphism is an open mapping (see \S
\ref{proof.Riemannan}). In particular
a harmonic morphism from a 1-dimensional objet (for instance a graph) to
${\Bbb C}$
is constant. The second part of the theorem above means that there are
non-constant holomorphic
functions on some graphs.



\subsection{Content of the article,  Rigidity in the discrete case}

Our first aim is to introduce this point of view of holomorphy, which we
find sufficiently motivated
as we said above, and it is at least amusing!
This point of view is obscure in the literature related to harmonic morphisms.
For example, (a)  Jacobi's problem, a precursor problem in the theory, since it
exactly
asks for a classification of global harmonic morphisms ${\Bbb R}^3 \to
{\Bbb R}^2$,
is expressed as a partial differential equation, but different from our
formulation of
holomorphic mappings (from ${\Bbb R}^3 \to {\Bbb C}$).  This problem was
solved only recently
\cite{B-W}. Other more ``elementary'' proofs were obtained afterwards (see for
instance \cite{Duh}). We anticipate
here to announce our ``elementary'' proof which will appear in a
forthcoming paper.

Also, as a contribution of the paper is the fact that the concept of
holomorphy in the discrete case seems to be
new. We were convinced of this for a long time, and  after parts of the
present paper were written,
we discovered (in special circumstances) that the terminology ``discrete
holomorphy''  was already used. 
Indeed in \cite{Mer} (and related papers) the author introduced a notion
of holomorphy for forms and functions on graphs.  Inspired by Hodge theory,
the definition is
global, and linear.  It is therefore different from our notion, yet,
possible analogies are
interesting to find.

 \subsubsection{Content.}   Instead of a systematic study of discrete 
holomorphy in the present
article, we
focus on  some particular examples which illustrate various  subtleties.


We will give some attention to the simplest non-trivial cases : 
the 3-valenced tree $T_3$
 and its (1-dimensional) ``dual graph'' $Tr_3$ (see \S \ref{TR3}).    
The spaces of holomorphic functions have
 finite dimension (which is of course not the case of harmonic
 functions),
and endowed with the beautiful action of the automorphism group.

 \subsubsection{Random dynamics versus Poisson kernel.}    In fact, 
holomorphic functions are ``essentially'' determined when
 given on some finite subsets of vertices.  There is a ``holomorphic dynamics'' which
 allows one to calculate the holomorphic functions progressively on the other
 vertices.  However, there is at each step some finite choice to consider,  in some
 sense, one has a kind of ``random dynamics''.

  One can conclude some general philosophy. In the case of harmonic functions,  
one goes
   to (space) infinity, and reproduce objects by a Poisson integral. 
In the holomorphic case,  a
   random holomorphic dynamics of finite subsets of the space itself 
allows one to reproduce
   the function.

  We do not claim our tentative to formulate all
 these notions is optimal. We find however exciting (and funny) 
the  beauty of these structures
 and their  connections with many natural notions, as correspondences (\cite{Bul}), ...



 Let us now summarize some obtained results.

 \subsubsection{The 3-valenced tree $T_3$.}

 The case of $T_3$ is the simplest non-trivial one.  Many (amusing) facts
 will be said about it.  The following theorem summarizes results obtained along
 \S \ref{T3}.

 \begin{theorem} \label{theorem.tr3}
 Let, $\phi \in Hol^*(T_3)$, that is,  a non-constant holomorphic function
 on $T_3$. Then:

 \begin{enumerate}

\item  $\phi$ sends (the vertices of) $T_3$ onto (the vertices of ) a hexagonal tiling
 $\tau$ of the Euclidean plane.

\item In fact
$ \phi$  is the universal covering
of the hexagonal tiling.  In particular:

$\bullet$ $\phi $ is locally injective, that is, it is injective on any (closed) ball
of radius 1 in $T_3$.
More precisely,
$\phi$   is conformal (see Remark \ref{conf} for definition).

$\bullet$  the image of a  geodesic of $T_3$ is a locally injective
walk on the hexagonal tiling. It could be periodic, but generally goes to infinity.


---  Let ${\mathcal H}_{\alpha,\beta}(T_3)$, for $\alpha \neq \beta$ be  
the space of holomorphic functions taking the values $\alpha$ and $\beta$ 
on two adjacent
vertices.  Let $H$ be the subgroup of $Aut(T_3)$ fixing (individually) these
vertices, then:


\item The action of $H$ on ${\mathcal H}_{\alpha,\beta}(T_3)$   
is simply transitive. In particular, the group and the space
are homeomorphic (in particular both are compact).


\item The same is true for
the product action of the group  $SG \times  H$ on the space $Hol^*(T_3)$,
where $SG$ is the similarity group of ${\Bbb C}$. 
The previous product group and space   are therefore homeomorphic.


\item  In contrast, the (transitive) action of $SG \times Aut (T_3)$ 
on $Hol^*(T_3)$ is not free. The stabilizer of a point 
(i.e. a non-constant holomorphic function) is a discrete group
$\Gamma$ acting freely on $T_3$, which is nothing but the fundamental group
of the hexagonal tiling.

 \end{enumerate}

  \end{theorem}

 One  sees from this how are rigid holomophic functions in comparison to
 harmonic ones. In particular, the space of holomorphic functions has finite
 (topological) dimension (exactly 4), whereas that of harmonic has
 infinite dimension.

  Also, one sees relationship between holomorphy and random walk on the hexagonal
  lattice, in fact, only  partial random walks are  involved, those,   called here locally injective.

 \subsubsection{The graph $Tr_3$.} Here, a holomorphic function is determined once
given on any triangle, but with some discrete indetermination, leading to a holomophic dynamics
on  the space of triangles.

 \medskip


$$\epsh{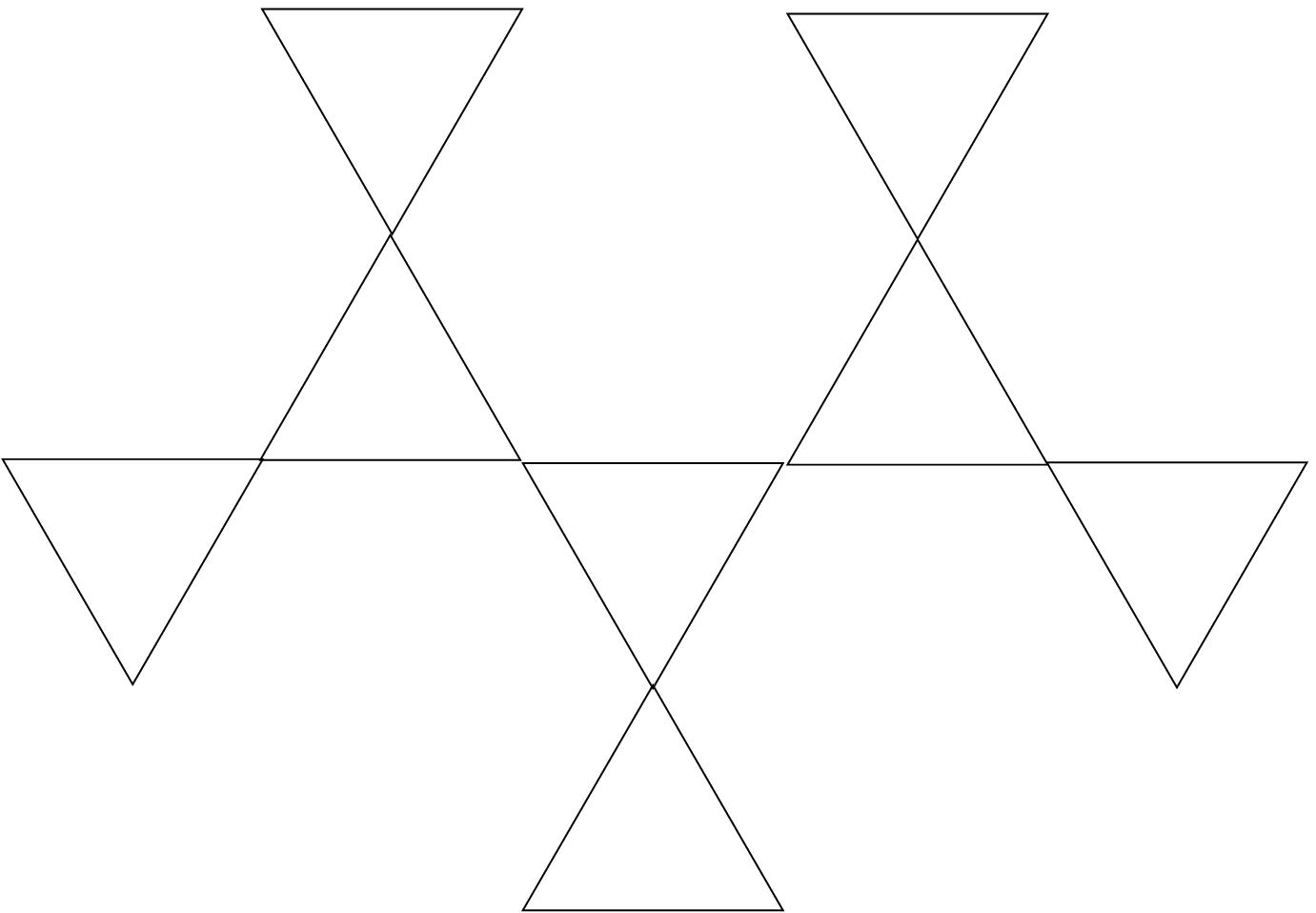}{3cm}$$ 
\begin{center}

{ Figure 1 :  The graph $Tr_3$}

 \end{center}

 \begin{theorem}  A holomorphic function on $Tr_3$  is determined once
given on any one of its triangles. There is however some  discrete indetermination
in the extending  process of the holomorphic function,  leading to a holomophic
random dynamics
in the space of triangles.

$\bullet$ {\em Correspondence.}  Precisely, 
this holomorphic dynamics is apprehended in the correspondence ${\mathcal S}$
 which is a 3-dimensional (complex) quadric defined in ${\Bbb C}^3 \times {\Bbb C}^3$
(endowed with coordinates $ (p, e, f; x, y, z) $)  by:
 \begin{eqnarray}
e+f-y+ ( -y+z) &=& 0 \nonumber \\
e^2+f^2+y^2+(-y+z)^2 &=& 0 \nonumber \\
p - (x+y) &=& 0 \nonumber
\end{eqnarray}

$\bullet$ {\em Random dynamics.} In other words, ${\mathcal S}$ is 
the graph of a holomorphic multi-valued transformation 
$M: (p, e, f) \in {\Bbb C}^3 \to (x, y, z) \in {\Bbb C}^3$, 
which encodes the diffusion process of a holomorphic function.

 $\bullet$ {\em Action of $Aut(Tr_3) $.} Two holomorphic 
functions are equal up to composition by an element
 of $Aut(Tr_3)$, iff, they take  same values on vertices of two triangles of $Tr_3$.

 $ \bullet$ {\em Orbital structure.}   
On the space of marked triangles ${\Bbb C} \times ({\Bbb C}^2-\{0\})$,  
consider the
 relation:  ${\bf \Delta} \sim {\bf \Delta^\prime} \iff \exists\;  \phi$ 
holomorphic
 such that $  {\bf \Delta^\prime} = \phi( {\bf \Delta})$. 
Then $\sim$ is an equivalence relation.

$ \bullet$ {\em Projectivization.}  By homogeneity,   we get a
correspondence in ${\Bbb C}P^1\times {\Bbb C}P^1 $, and an equivalence relation in
${\Bbb C}P^1$.
 \end{theorem}

(It turns out that the equivalence relation $\sim$, 
or its projective variant, are not defined by group actions).

These two theorems do not answer
many other remaining  questions on the subject of holomorphic functions  
on   these special two spaces  $T_3$ and $Tr_3$. Actually, 
in this introductory article, we could
not go deeply in analysis : the phenomenology by itself is rich enough, 
and formulations are not straightforward!  Preliminary 
simulations promise a beautiful and strong dynamics.

 \subsubsection{Other graphs, non-rigid situations.} The results for $T_3$ allow one
to understand holomorphic functions for any 3-valenced graph. This is discussed
in \S \ref{remarks.other}, where also  remarks are made, 
giving evidences to a rigidity
phenomenon in the case of the Cayley graph of ${\Bbb Z}^2$.  In contrast, one
observes abundance of holomorphic functions for trees of valency $\geq 4$.
The study of these trees will be continued in \S \ref{conjugate.problem}, where
we formulate the ``conjugate part problem''.

\section{Some properties of harmonic morphisms, Proofs in the Riemannian case}
\label{proof.Riemannan}

We recall in what  follows,  fundamental properties  of harmonic morphisms
in the context of
Riemannian manifolds. They are due   to  B. Fuglede and T.  Ishihara
(independently) \cite{Fug, Ish}. Let
 $\phi: X \to Y$ be a smooth mapping between two Riemannian manifolds. Then
{ \it $\phi$ is a harmonic morphism, iff, $\phi$ is harmonic and horizontally
semi-conformal. } The latter means that, 
if $\phi$ is not constant, then 
for almost all  $x$, if $H_x$
denotes the normal
space of $\ker D_x\phi$, then: $D_x\phi: H_x \to T_{\phi(x)}Y$ is a
conformal isomorphism
(on the complementary negligible set where $D_x\phi = 0$).

It then follows in particular that a harmonic morphism is an open mapping.

As an example if $\dim X =1$, and $\dim Y = 2$, then the harmonic morphism
$\phi$ is constant.

 \begin{remark}
 {\em
This fact can also be proved (in an essentially similar  way) when $X$ is merely a
graph (say, a 1-dimensional object
more general than manifolds). This explains the comment after Theorem
\ref{non.quadratic}.
 }

 \end{remark}

Consider now a mapping $\phi: X \to {\Bbb C}$, where $X$ is a Riemannian
manifold.

Assume $\phi$ is a harmonic morphism. The mappings $z \to z^N$ are harmonic
on ${\Bbb C}$. By definition,
$\phi^N$ is harmonic for any $N$, in particular $\phi$ is holomorphic.

Assume now that $\phi$ is holomorphic.  From the discussion following
Definition \ref{definition.holomorphic},
it follows that $\phi$ is horizontally semi-conformal (as defined above).  
>From the above
characterization,
$\phi$ is a harmonic morphism.

Summarizing, $\phi$ is a harmonic morphism, iff,
$\Delta \phi = \Delta \phi^2=0$, iff, $\Delta \phi^N =0$ for all $N$.

This proves the Riemannian parts in Theorems \ref{quadratic} and
\ref{non.quadratic}

Finally, let us say that the fact that harmonic morphisms satisfy a
quadratic equation,
 or equivalently that $(*)_2 \Longrightarrow (*)_\infty$, is stated through
standard papers like \cite{Fug}.
However, it was just  used as an
intermediate step of
proofs.

 \section{Generalities on graphs}

 We start with a  general local  study,  and fix  notations.
Let $s_0 $ be  a vertex of $X$, with valency $n$: $s_1, \ldots, s_n$ 
are the adjacent vertices of $s_0$.
The formel difference $s_i-s_0$ denotes the oriented edge joigning $s_0$ to $s_i$.

Let $\phi$ be a (complex) function on $X$, with $z_i= \phi(s_i)$. 
The oscillation of $\phi$ along
the  edge $s_i-s_0$ is $ \delta_i= z_i-z_0$.
The gradiant of $\phi$ at $s_0$ is  the vector 
$\overrightarrow{\delta}=(z_1-z_0,\ldots,z_n-z_0) \in {\Bbb R}^n$.
(Actually the values $\{z_i\}$ are not ordered and therefore, 
$\overrightarrow{\delta}$ is an element of the symmetric
quotient of ${\Bbb R}^n$. However this precision does 
not matter in our analysis here, and for the sake of
simplicity, we keep our notations above).

The function $\phi$ is harmonic at $s_0$ if its gradiant at $s_0$ 
has a vanishing divergence (i.e.
arithmetic mean): $\sum(z_i-z_0) = \sum \delta_i= 0$.

 \begin{fact} \label{oscillation}  The function $\phi$ is holomorphic,  iff,  both its mean oscillation and
 quadratic  mean oscillation vanish at any vertex:
 $$\sum \delta_i = \sum \delta_i^2 =0$$
 In fact, more generally, for a harmonic function $\phi$,  the mean quadratic oscillation equals
 the mean oscillation of its square $\phi^2$.

 \end{fact}

 \begin{proof}   At $s_0$, the mean oscillation of $\phi$ and $\phi^2$ are respectively:
 $ \sum \delta_i= \sum (z_i-z_0)$ and $\sum(z_i^2-z_0^2)$.

  Assume $\phi$ is harmonic at $s_0$, then,
 $\sum z_i = nz_0$ ($n$ is the valency at $s_0$).  Thus,
 $$\sum \delta_i^2 = \sum(z_i-z_0)^2= \sum z_i^2+ nz_0^2-2z_0\sum z_i =
 \sum z_i^2+nz_0^2-2nz_0^2 = \sum(z_i^2-z_0^2)$$

 \end{proof}

 We will also need the following lemma. Its proof reduces to the possibility of solving
  in ${\Bbb C}$,
 a system of a linear and a quadratic equations, on two unkowns.

 \begin{lemma} \label{lemme}
Let $n\geq 3$ and  $\delta_1,...,\delta_{n-2}\in \Bbb C$ be given.   Then there
exist a  pair $\lbrace \delta_{n-1} , \delta_n \rbrace$,
 unique  up to switch,
  of (complex) solutions of
$$(\sum_{i=1}^{n-2} \delta_i)+\delta_{n-1} + \delta_n =
(\sum_{i=1}^{n-2} \delta_i^2)+\delta_{n-1}^2+ \delta_n^2= 0.$$

 \end{lemma}

 In particular, to fix notations for future use in \S \S \ref{T3} and  \ref{TR3}:

 \begin{fact}  \label{solution}
 For given $e$ (resp. $e$ and $f$)
 there is  $\lbrace u,v  \rbrace$  a  unique pair up to switch, such that:

\begin{eqnarray}
e+u+v &=&0   \nonumber \\
  e^2+u^2+v^2&=&0 \nonumber
\end{eqnarray}


( resp.
\begin{eqnarray}
e+f+u+v &=&0   \nonumber \\
  e^2+f^2+u^2+v^2&=&0 \nonumber \;) \label{eq:equa}
\end{eqnarray}

\end{fact}

  \section{Dynamical and ergodic study for the 3-valenced tree $ T_3$} \label{T3}


 We investigate here holomorphic functions on $T_3$, the (bi-infinite) tree of valency 3. However, because of technical difficulties (due to the fact that $T_3$ can not be well ordered, that is endowed
 with an orientation inducing an orientation on any  geodesic), we will start by considering
 intermediate cases.

 \subsection{The smallest tree, the tripod}

This is a  (finite) tree ${\bf Y}$ with 4 vertices, $O, O^\prime, A,B$, all of valency 1, except $O$, which
has valency 3 (all the others are extremal).

  Let $\phi$ be a function on ${\bf Y}$ , and denote by  $e, u, v$ its  associated oscillations at $O$, that is
 $e= \phi(O^\prime)-\phi(O)$, $u= \phi(A)-\phi(O)$ and $v = \phi(B)- \phi(O)$.

 Suppose $\phi$ is holomorphic at $O$, then, $$e+u+v= e^2+u^2+v^2=0.$$

 Fact \ref{solution} says that $u$ and $v$ are completely  determined, up to order,  when $e$ is given. More precisely, let
 $j = e^{\frac{2}{3} \pi i}$ be  a cubic root of unity, then, necessarily, up to order, $$u= j e, \quad v = j^2 e$$

 For instance,
 $\phi(A)- \phi(O)=  \alpha (\phi(O)-\phi(O^\prime)) $
  for some $\alpha \in \{-j, -j^2\}$.

  Thus,  $$\phi(A)=  \phi(O) + \alpha (\phi(O)-\phi(O^\prime))$$

Also, $\phi(A)= \phi(O^\prime)+ (\phi(O)-\phi(O^\prime)) + 
\alpha (\phi(O)-\phi(O^\prime)$.


 \medskip

$$\epsh{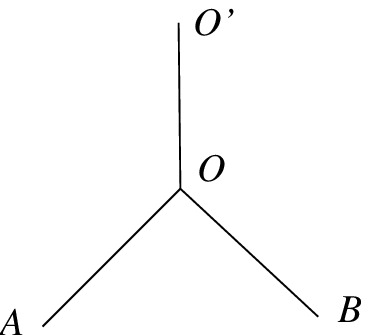}{2.5cm}$$ 



  \begin{center}
  { Figure 2 : The tripod}
 \end{center}


 \begin{remark} \label{conf} (Conformality)  {\em
  The natural ``geometric'' way to embed ${\bf Y}$
 in the Euclidean plane is to imagine the edges $OO^\prime, OA$ and $OB$,  having
 the same length and making equal angles at $O$, which must be $\frac{2 \pi}{3}$.  We see from above,
 that a holomorphic mapping $\phi$ is then conformal (the target space ${\Bbb C}$ is of course endowed with its
 Euclidean structure).

 }

 \end{remark}

 \subsection{Rooted  tree}  \label{marked.tree}

 Let now ${\bf A}$ be an infinite tree of valency 3, with a root $R$. So, all vertices have valency 3, except
 $R$ which has valency 2.   A model of it is given by
 the free monoid on two letters $\{a,  b\}$, the root $R$ corresponds to the void word. The vertices are words on $a$ and $b$.

 Instead of ${\bf A}$, we consider an extended tree ${\bf A}_{O^\prime O}$, obtained by gluing
 an  oriented segment $O^\prime O$ to ${\bf A}$, where $O$ is identified with $R$. Now, the root becomes
 $O^\prime$.

 If $\phi$ is a function given  at  $O^\prime$ and $O$, then one extends it holomorphicaly step
by step, applying the rule described for ${\bf Y}$. To formulate this, orient
${\bf A}_{O^\prime O}$ naturally, in such a way that $O^\prime$ is a source, and for every
$S$ a word on $\{a, b\}$,  the edges $Sa$ and $Sb$ are positively oriented.

If, $\phi$ is already  defined on $S^\prime$ and $S$, so that $S^\prime S$
is an oriented edge, then,
$$\phi(SX)= \phi(S) + \alpha_S(X) (\phi(S)- \phi(S^\prime)),  \;  X \in \{a, b\}$$
 where $$\alpha_S: \{a, b\} \to  \{-j, -j^2\}$$
   is a {\it bijection}.

 Therefore, $\phi$ is completely determined by, the ``randonm variables'' $\alpha_S$, associated
 to any $S$ (a word on $a$ and  $b$).

$$\epsh{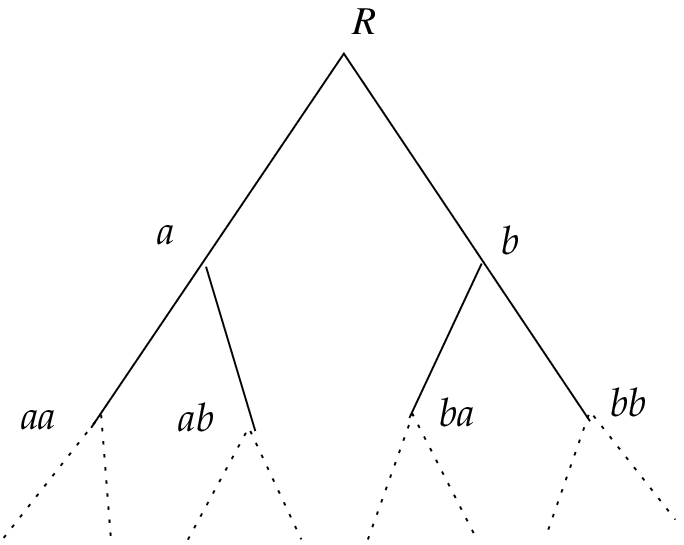}{4cm}$$ 
\begin{center}

 { Figure 3 : Rooted tree}

 \end{center}

 \subsubsection{Chain rule} \label{chain}

If $O^\prime= S_{-1},  O = S_0,\ldots S_n$ is a geodesic, i.e. a simple path   in ${\bf A}_{O^\prime O}$, and
 $w = \phi(O)-\phi(O^\prime)$,  then,
 $$ \phi(S_n) = \phi(O^\prime) +  w+  
\alpha_0w+ (\alpha_1\alpha_0) w +\ldots+
( \alpha_{n-1}\alpha_{n-2}. \ldots . \alpha_0)w,$$
  where $\alpha_i  = \alpha_{S_i} $.

 \subsubsection{Planar orientation, Canonical example}

 For two given different complex numbers $\alpha$ and $\beta$,  there is a  canonical holomorphic function  $\Phi$ with $\Phi(O^\prime) = \alpha$ and $\Phi(O) = \beta$.
 It is   determined by taking a same $\alpha_S$ for all $S$, defined by:
 $$\alpha_S(a)= -j, \; \alpha_S(b)= -j^2$$

This $\Phi$ can be characterized among holomorphic mappings by  the fact that it preserves
``2-dimensional'' (or planar) orientation of the tree.  To define it, imagine the tree naturally embeded
in an oriented Euclidean plane. 
This induces an ``orientation'' on the tree, that is a choicee of positive
angle between $Sa$, and $Sb$, for any $S \in {\bf A}$. Here, we decree that
$\angle(Sa, Sb)$ is positive.  In this case $\Phi$ actually preserves the orientation. It will be seen just
below how to obtain all holomorphic mappings by means of $\Phi$, and then that only $\Phi$ preserves orientation.

 Summarizing these nice properties of $\Phi$:

 \begin{fact}  \label{nice}  $\Phi$ is holomorphic, (planar) orientation preserving,
 and {\em locally injective}, that is,  at any vertex,
 $\Phi$ is injective on its immediate neighbourhood (i.e. on the ${\bf Y}$ around it), and is conformal
 (as defined in Remark \ref{conf}).

 \end{fact}

  \subsubsection{Simply-transitive  action of  $Aut$}
   Let ${\mathcal H}_{\alpha,\beta}$ be the space of holomorphic functions on
   ${\bf A}_{O^\prime O}$ with prescribed  (and distinct) values  $\phi(O^\prime) =\alpha$ and
    $\phi (O) = \beta $.

    Let $Aut ({\bf A})$ be the automorphism group of ${\bf A}$ or equivalently   
 ${\bf A}_{O^\prime O}$.
    It fixes $O^\prime$ and $O$, and for any $S$, there 
are elements of the group which
    exchange $Sa$ and $Sb$.
     On the other hand, this exchange was the only one 
freedom in determining the holomorphic function, step by step, as explained above.
    This can be translated into the fact that, the right action by composition of
         $Aut({\bf A}_{O^\prime O})$   on ${\mathcal H}_{\alpha,\beta}$, is transitive.

  On the other hand, this action is free.  
To see this, one verifies that $\Phi$ has a trivial stabilizer, that is no
  $\gamma$ except the identity   satisfies : 
$\Phi \circ \gamma = \Phi$.  To prove it, one sees that, such a
  $\gamma$ must  fix  both $a$ and $b$. Therefore, at the next step, 
that is for words of length 2, the only
  non-trivial possible action of $\gamma$ is to induce 
a transposition  on $\{aa, ab\}$
  or $\{ba, bb\}$.  This means $\Phi$ takes a unique value on 
one on these sets, which contradicts
  the fact that $\Phi$ is locally injective (Fact \ref{nice}). 
Now, the formal proof of the freedom of the action can
  be performed by induction.

Both  ${\mathcal H}_{\alpha,\beta}$ and $Aut({\bf A}_{O^\prime O})$ are compact
(when endowed with the compact-open topology,
 although the elements of ${\mathcal H}_{\alpha,\beta}$ are not bounded!).
They are in  fact topologically
      Cantor sets.
      One checks that: $\gamma \to \Phi \circ \gamma$ is a homoeomorphism.

  \subsubsection{Action  of the similarity group} 
The similarity group $ SG$  of ${\Bbb C}$ acts (by the left)
  on the space
  of holomorphic functions for any space. They are mappings : 
$z \to \theta z+t$, where $\theta, t \in {\Bbb C}$,
  $\theta \neq 0$.

  In our case,  let $Hol= Hol({\bf A}_{O^\prime O})$ be 
the space of holomophic functions. This is a locally compact
  space (for the compact-open topology).  Let $Hol^*$ be the set of non-constant functions, that is,
  $Hol^*= Hol - {\Bbb C}$. Equivalently, $Hol^*$ is the (disjoint) union of all
  the ${\mathcal H}_{\alpha,\beta}$, for $\alpha, \beta \in {\Bbb C}$, $ \alpha \neq \beta$. Then, $SG$ permutes the spaces ${\mathcal H}_{\alpha,\beta}$. In fact, a non-trivial
  element of $SG$ can  preserve no individual ${\mathcal H}_{\alpha,\beta}$.  On the other hand, $Aut({\bf A}_{O^\prime O})$ acts by preserving
  each ${\mathcal H}_{\alpha,\beta}$. One then consider the produt action of  $SG \times Aut({\bf A}_{O^\prime O})$. From
  previous developments, we get:

  \begin{fact} The product action  of $SG \times Aut({\bf A}_{O^\prime O})$ on $Hol^*$ is simply transitive. Consequently,   $Hol^*$ is homeomorphic to $SG \times
  Aut({\bf A}_{O^\prime O})$.

   \end{fact}


 \subsection{Full tree $T_3$}

The full (bi-infinite) tree $T_3$ is obtained by gluing  ${\bf A}_{O^\prime O}$ and ${\bf A}_{O O^\prime}$ along their respective edges
$O^\prime O$ and $O O^\prime$.  The group $ Aut({\bf A}) \times Aut({\bf A})$ can be viewed as
the subgroup of $Aut (T_3)$ fixing the edge $O^\prime O$ (i.e. each of its extremities). With a little bit analysis, one gets:

\begin{proposition}  \label{action.full}
 An element of $Hol^*(T_3)$, that is a non-constant holomorphic function 
on $T_3$ is locally injective,
more precisely it is conformal (as defined in Remark \ref{conf}).

Let ${\mathcal H}_{\alpha,\beta}(T_3)$, for $\alpha \neq \beta$, 
be the space of holomorphic functions taking values $\alpha$ and $\beta$ 
on two adjacent
vertices. The action of  $ Aut({\bf A}) \times Aut({\bf A})$ on it is 
simply transitive. The same is true for
the action of the group  $SG \times  Aut({\bf A}) \times Aut({\bf A})$ 
on the space $Hol^*(T_3)$, which are
therefore homeomorphic.

In contrast, the (transitive) action of $Aut (T_3)$ on $Hol^*(T_3)$ 
is not free. The stabilizer of a point is a
(canonical) discrete group
$\Gamma$ acting freely on the tree $T_3$.

\end{proposition}



 \subsection{Hexagonal tiling}

To given  $\alpha$  and $\beta= \alpha+ w$, one associates a hexagonal tiling of
 ${\Bbb R}^2 $ (=${\Bbb C}$), with $\alpha$ as a vertex, and $w$ as an edge (trough $\alpha$). The fundamental
 tile is the regular hexagon with vertices 
$\alpha, \alpha+w, \alpha+w-jw, \alpha+w-jw+j^2 w, \alpha+w-jw+j^2 w-w,
 \alpha+w-jw+j^2 w-w+jw$.

 The tiling $\tau_{\alpha,\beta}$ is obtained by reflections along its edges of the fundamental
 hexagon.  Up to similarity, we can restrict ourselves to $\tau= \tau_{0,1}$.

 From \S \ref{chain}, we get:

 \begin{proposition} Let, $\phi \in Hol^*(T_3)$, that is a non-constant holomorphic function
 on $T_3$. Then, up to rescaling by a similarity of ${\Bbb C}$, we have:

$ \bullet$ $\phi$ sends (the vertices of) $T_3$ onto (the vertices of ) $\tau$.

$\bullet$  In particular, for a geodesic (injective path)
 ${\mathcal C}= S_0, S_1, \ldots, S_n, \ldots$ in $T_3$,  the image, $c = s_0,\ldots, s_n= \phi(S_n), \ldots$ is a path of $\tau$.

- The path $c$ is locally injective, that is $s_i \neq s_{i+1} \neq s_{i+2}$.

 $\bullet$ Moreover, because of the local injectivity of $\phi$, a ramification of paths in $T_3$ induces a
 ramification in $\tau$. More exactly, if two geodesics ${\mathcal C}= S_1,\ldots, S_n,  \ldots$ and
 ${\mathcal C}^\prime= S^\prime_1, \ldots, S^\prime_n, \ldots $ bifurcate at $i$,  that is,
 $S_j= S^\prime_j$, for $j<i$, and $S_i \neq S^\prime_i$, then, the same is true for their
 images, that is,  $s_i \neq s^\prime$.

 $\bullet$ Summarizing:  $\phi: T_3 \to \tau$ is the universal covering, and  the discrete  group
 $\Gamma$ in Proposition \ref{action.full} is the fundamental group of
 $\tau$

 \end{proposition}


   \medskip

$$\epsh{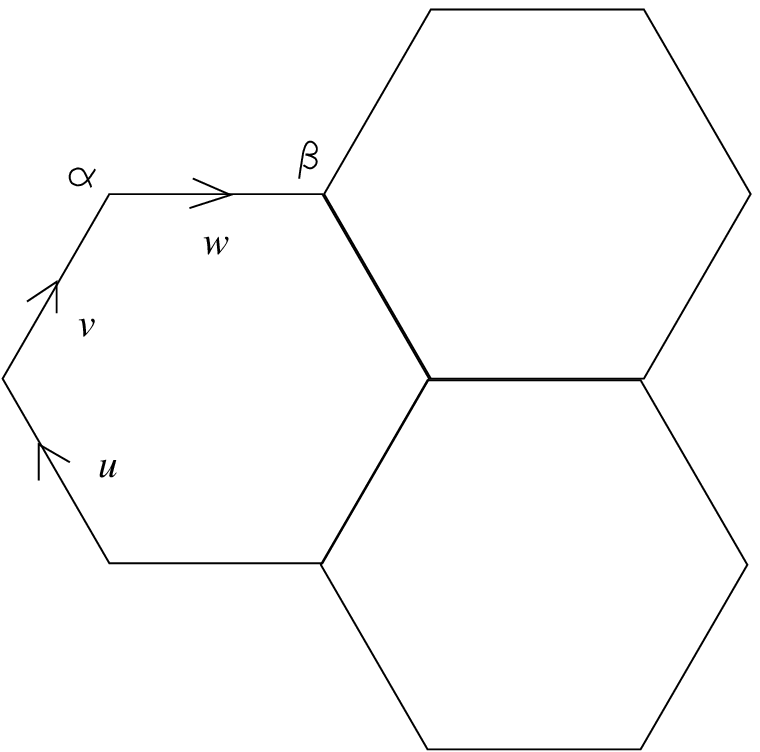}{4cm}$$ 
\begin{center}

{ Figure 4 : Hexagonal tiling}

 \end{center}

   \subsubsection{Locally injective random walk on $\tau$}

   Let ${\mathcal A} $ be the set of  oriented edges of the fundamental hexagon of $\tau$,
   which are labelled $u, v, w$, and then
   (necessarily) $-u, -v, -w$.
   Consider $\Omega= {\mathcal A}^{\Bbb N}$. Let $X_i:  \omega = (\omega_i)
   \in \Omega \to \omega_i \in {\Bbb C}$.
    Consider (partial) sums
  $S_i(\omega)= \sum_{j=0}^{j=i} X_j(\omega)$. A 
sequence $S_0, S_1, \ldots S_n$ is supposed to represent a path of $\tau$. 
However, not all $\omega$'s are allowed. For instance, if $S_0(\omega)= u$, 
then, in order
  that $\omega$ represent a walk in $\tau$, $X_1(\omega)$ must equals $ -u, v$ or $-w$.
  In general $\omega$ belongs to a subshift of finite type defined by means of an incidence matrix $A^\tau$, with entries 0 and 1: a sequence $\omega $ is allowed, iff, $A^\tau_{\omega_i\omega_{i+1}}
  = 1$, for all $i \in {\Bbb N}$.  Now, being furthermore a {\it locally injective} path in $\tau$, leads
  to  another subshift $\Sigma^\tau_{inj}$. The new forbidden subsequences are of the form
  $y, -y$, where, $y \in {\mathcal A}$. If ${\mathcal A}= \{u, v, w, -u, -v, -w\}$ is identified, preserving order with
  $\{1, \ldots, 6\}$, then the incidence matrix of $\Sigma^\tau_{inj}$ is:

\[
\left(
\begin{array}{cccccc}
 0 & 1& 0& 0 &0 &1 \\
 1&0&1&0&0&0 \\
  0 & 1 & 0 & 1 & 0 & 0 \\
    0 & 0&1&0&1&0   \\
  0&0&0&1&0&1  \\
 1 & 0 & 0 & 0 & 1 & 0
\end{array}
\right)
\]


 $\bullet$ Ergodic theory of the subshift $\Sigma^\tau_{inj}$, growth of random variables $S_n$, are
 straightforwardly   related to growth of holomorphic functions on $T_3$.













\section{Proof of the graph case in Theorems \ref{quadratic} and \ref{non.quadratic}}

\label{graph.case}

We have the following straightforward generalization of Lemma \ref{lemme}

\begin{fact}

Let $\delta_1, \delta_2$ be given. Consider the set of equations on unknowns
  $\delta_3,...,\delta_n$:

$$\sum_{i=3}^{n} \delta_i^p+ (\delta_{1}^p + \delta_2^p) = 0, $$
for $p \in \{1, \ldots, N \}$, where $N$ is an integer.

If $N \leq n-2$, then a solution exists. 
It is unique up to permutation (of the $\delta_i, i\geq3$)
if $N = n-2$. In opposite, if $N>n-2$, then, 
in order that solutions exist, $\delta_1$ and
$\delta_2$ must vanish, and the solution in this case is trivial : 
$\delta_i =0,\; \forall i$.
 \end{fact}

This yields as a corollary  the claim  that a holomorphic function which
satisfies $(*)_N$  where $N$  is the (maximal) valency
of the graph $X$  is constant.  In particular, in the case of $T_3$ only constant
functions satisfy  $(*)_3$, so all other holomorphic functions satisfy $(*)_2$ but not
$(*)_3$.

A harmonic morphism $X \to {\Bbb C}$, must satisfy $(*)_\infty$. Only constant functions are
so (say, assuming the graph has a finite valency).

This proves all the discrete content of Theorems \ref{quadratic} and \ref{non.quadratic}.

 \subsubsection{$N$-holomorphic functions.}  The previous fact suggests to introduce 
$N$-holomorphic
 functions as those satisfying $(*)_N$.  For instance a theory of 
$(*)_{N}$ holomorphic functions on
 $T_{N+1}$, the tree of valency $N+1$ may be developed in a same way as 
holomorphic functions on
 $T_3$. The situation is  a  little bit complicated. For instance, 
instead of the hexagonal tiling, we get
  a picture with other kinds of piecewise lines. (See the following figures 
in the cases $N=2,3,4,5$)

$$\epsh {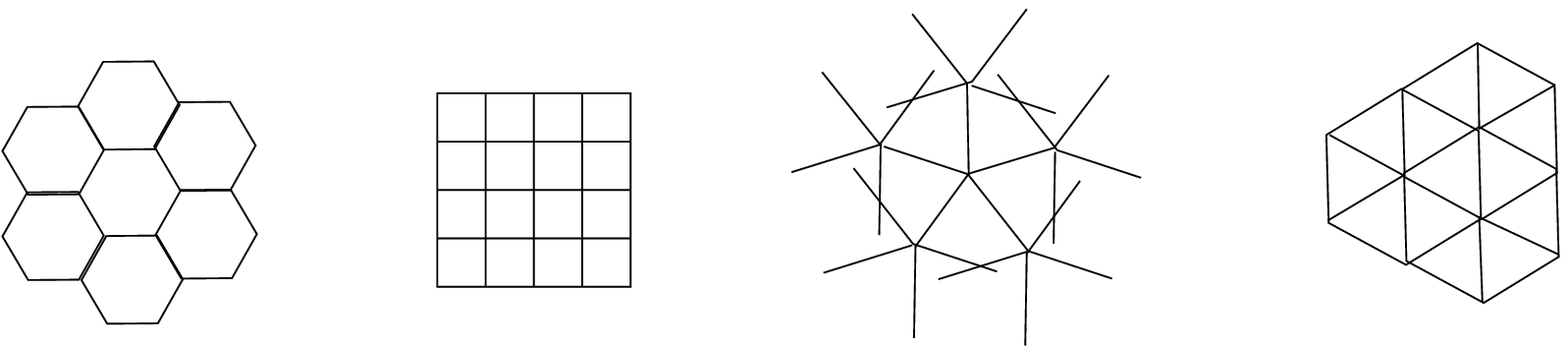}{2.5cm}$$ 
\begin{center}

{ Figure 5 : N-holomorphic functions }

 \end{center}


\section{The graph $Tr_3$} \label{TR3}

The graph $Tr_3$, is a  kind of 3-valenced tree, where vertices are replaced by triangles, which
explains our notation here.

More precisely $Tr_3$ is a homogeneous graph with valency 4 and (minimal) cycles of length 3.

It is obtained from a countable family of triangles, where, to  each triangle is glued exactly another
one at each vertex.

It can also be considered as the ``dual graph'' of  the tree $T_3$, that is,  the graph which vertices are the edges of $T_3$, and two of them are joined by an edge (in $Tr_3$) if the have a common vertex (in $T_3$).  A concrete model is obtained by taking middles of the edges of $T_3$ and drawing segments joining middles
of adjacent edges.


In the sequel, we will see, analogously to the case of $T_3$, that a holomorphic 
function, is constructed
following some dynamics, when it is given at vertices of some triangle. 
As  in the case of $T_3$,
a nuisance  but  also a richness of the dynamics comes from its 
non-deterministic character,
 that is,  the switch choice.    However,
with respect to this point, the case of $Tr_3$ is highly  
more complicated. Our attempt here is just
to roughly study  this dynamics, which we think is worthwhile 
to investigate, at least  because
it is natural!

\subsection{Dynamics on triangles}

\subsubsection{Starting-up}
Let $OAB$ be  a triangle in $Tr_3$, and $OCD$ and adjacent one, 
with a common vertex $O$.  If $\phi$ is holomorphic, and is given  
at $O, A$ and $B$, then  using
Fact \ref{solution} in a same way as in the case of $T_3$, one deduces
the values of $\phi$ at $C$ and $D$, up to a switch.

Let's introduce notations which will be kept along this \S. We denote:
$$p= \phi (O),
e= \phi(A)-\phi(O),  f= \phi(B)-\phi(O),
u= \phi (C)-\phi (O), v=  \phi(D)-\phi(O)$$

Here $p, e$ and $f$ are given. We infer from Fact \ref{solution}, that $u$ and $v$ can be computed, but up to
a switch. Indeed, this reduces to solving the following equation on $u$ and $v$, where $e$ and $f$ are supposed known:
$$u+v+ (e+f)= u^2+v^2+ (e^2+f^2)= 0$$

In other words,  after elementary algebraic manipulation:

  \begin{fact} \label{second}
  $u$ (or equivalently $v$) is a solution of the second order equation:
  $$ u^2+ (e+f)u+ (e^2+f^2+ef)=0$$
   {\em (In the sequel, we will choose, arbitrary  solutions,   $\psi_1(e, f)$ and $\psi_2(e, f)$, see below)}.

  \end{fact}


  \medskip
$$\epsh{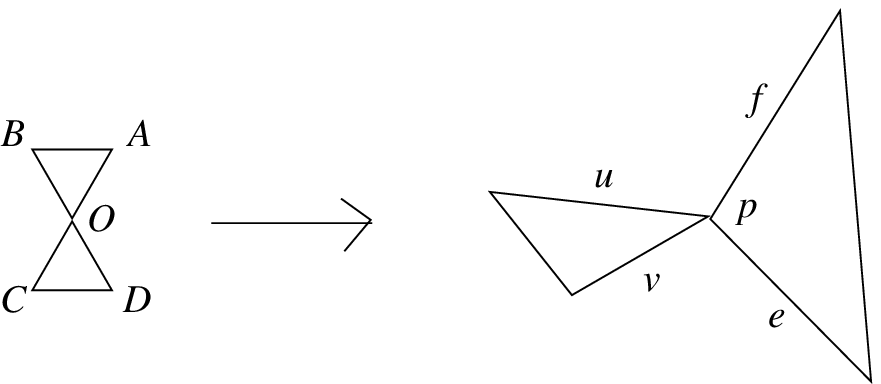}{3cm}$$ 

\begin{center}
{ Figure 6 :  Two triangles}

 \end{center}

\subsubsection{Rough formulation} \label{rough}

Now, once $\phi $ is determined on $OCD$ one uses data on this triangle to extend $\phi$ on
 an adjacent triangle, say that with vertex $C$. The data are thus
 $(\phi(C), \phi(O)-\phi(C), \phi(D)-\phi(C))$, or equivalently $(p-u, -u, -u+v)$.  Therefore, we are led
 to study the transformation:
 $$M:  (p, e, f) \to (p-u, -u, -u+v)$$

 The difficulty comes of course from that this is not  an  univalued transformation.
 Indeed $M$ is constructed by means of  the 2-valued mapping:
 $$I: (e, f) \in {\Bbb C}^2 \to (u, v) \in {\Bbb C}^2$$

In order to get univalued mappings, we choose, in an arbitrary way,  solutions   $\psi_1, \psi_2: {\Bbb C}^2 \to {\Bbb C}$
 as in  Fact \ref{second}.
 Actually, at this stage, we do not mind on any regularity of these mappings., but the only ``algebraic'' condition is that,  this gives rise to an involution, that is,  $\psi =(\psi_1, \psi_2)$ is an involution.

 Therefore, we get
 two ``sections'' for $M$:
 $$M_1: (p, e, f) \to (p- \psi_1(e, f), - \psi_1(e, f), - \psi_1(e, f)+ \psi_2(e, f)), $$
 $$M_2: (p, e, f) \to (p- \psi_2(e, f), - \psi_2(e, f), - \psi_2(e, f)+ \psi_1(e, f))$$

Summarising:

\begin{fact} \label{triangle.dynamics} If the image of the triangle $(O, OA, OB)$ by $\phi$
is $(p, e, f)$, then the image of  $(C,  CO, CD)$  and $(D, DO, DC)$ are, in order,
either,
$M_1(p, e, f)$, $M_2(p, e, f)$, or inversely $M_2(p, e, f)$, $M_1(p, e, f)$.
 \end{fact}

\subsubsection{The involution}

An important part of dynamics is  contained in
 the 2-valued mapping,
$I: (e, f) \in {\Bbb C}^2 \to (u, v) \in {\Bbb C}^2$.
It is an involution (as a 2-valued mapping).

One way  to let this being a ``standard'' mapping, is to consider the quotient space
${\Bbb C}^2/{\Bbb Z}_2$, where ${\Bbb Z}_2 = {\Bbb Z}/2{\Bbb Z}$ acts by $(e, f) \to (f, e)$.
This gives rise to a well defined (continuous) involution $\overline{I}$
of ${\Bbb C}^2/{\Bbb  Z}_2$.

In order to see a simple picture of this, it is useful to projectivize all things.  In this case,
the ${\Bbb Z}_2$ action on ${\Bbb C}P^1$
 becomes a $\pi$ rotation, around, north and south poles, say. The quotient space is still, topologically,
 a 2-sphere.

 It is suggestive to think about $\overline{I}$ as a holomorphic mapping of the ``orbifold''
 ${\Bbb C}^2/{\Bbb Z}_2$ (or equivalently ${\Bbb C}P^1/{\Bbb Z}_2$). This seems to not be so right, since
 for instance, $\overline{I}$ does not preserve the singular locus of orbifolds.
Anyway, the true dynamics which interests us, that is, $M_1$ and $M_2$, are not compatible
with tacking this quotient.

Maybe, the induced mapping ${\Bbb C}^2 \to {\Bbb C}^2/{\Bbb Z}_2$,  could be thought, 
more naturally,
as holomorphic.  However, having a different source and target spaces,  
deprives one to perform
dynamics (i.e. to iterate).

 \subsubsection{Catastrophe.}  The ``singular locus'' of $I$ is the set of  $(e, f)$, where
equation \ref{solution} has a double root, that is $u= v$. 
It  consists of the two complex lines
${\mathcal S}$ and ${\mathcal N}$ in ${\Bbb C}^2$
generated by $(1-i\sqrt{2}, 1+i\sqrt{2})$ and $(1+i\sqrt{2}, 1-i\sqrt{2})$. In the projective plane
 ${\Bbb C}P^1$, this  corresponds to two points, $\{S, N\}$, say. 
One then hopes that
 $\psi_1$ and $\psi_2$ (introduced in \S \ref{rough}) 
can be continuously defined outside
 ${\mathcal S} \cup {\mathcal N}$, or equivalently on $ {\Bbb C}P^1-\{S, N\}$.  
It turns out that this is not possible : despite $\psi_1 \neq \psi_2$ 
along $ {\Bbb C}P^1-\{S, N\}$, the {\it real}  line 
bundle that they define (the fiber above $(e, f)$ is the real line joining
 $\psi_1(e, f)$ and $\psi_2(e,f)$)
 is not (topologically) trivial, i.e. non-orientable (on $ {\Bbb C}P^1-\{S, N\}$).

\subsection{Random dynamics}

The graph $Tr_3$ can be encoded as follows. We have a central triangle $\Delta_0$.  
We denote by $1, 2, 3$ its vertices.
Any triangle $\Delta$ in $Tr_3$ is obtained by specifying a 
vertex $i$ which indicated the ``direction'' of $\Delta$, together
with a word of  the free monoid $F_2$ on two letters $\{a, b\}$ 
(which may be empty $\emptyset$). Therefore, any $\Delta \neq \Delta_0$ 
is encoded by
$(l, i)$, where $l \in F_2$, and $i \in \{1, 2, 3\}$.

For instance, in the direction 1, we have: 
$$\Delta_0, (\emptyset, 1), ( (a, 1), (b, 1)), 
( (aa, 1), (ba, 1), (ab, 1), (bb, 1)), \ldots$$


$$\epsh{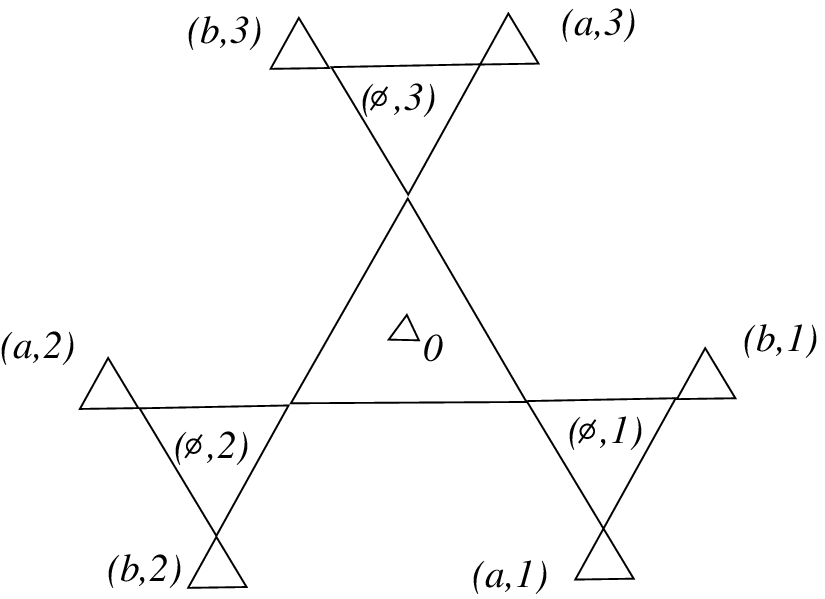}{4cm}$$ 

\begin{center}

{ Figure 7 }

\end{center}

Actually, one must  deal with marked triangles. For instance, the central (geometric) triangle $\Delta_0$ gives rise
to three marked triangles $(\Delta_0, i)$, where $i \in \{1, 2, 3\}$ indicates the vertex.

Each triangle $(l, i)$ gives rise to 3 marked triangles, one for each vertex.  Firstly,  a marked triangle $(l, i, 0)$, where the marked vertex is that sheared
by its (previous) generating triangle.
Now, imagine $Tr_3$ naturally embed in an oriented Euclidean plane. This induces an orientation
(i.e. an order on the edges) of the marked triangle $(l, i, 0)$. We get therefore two marked triangles
$(l, i, -)$ and $(l, i, +)$.

A generic marked triangle will be denoted $\overline{\Delta}$ and its underlying (non-marked) triangle is denoted
$\Delta$.

One can speak of $c \Delta$, for
$c $ a letter ($\in \{a, b\}$).
Also,  one defines, in a natural way, $c\overline{\Delta} = c (l, i, \epsilon)= (cl, i, \epsilon)$, for $\epsilon \in \{0, -, +\}$ a sign.
In other words,   the action of $F_2$ is trivial on the marking.





 Exactly as in \S \ref{marked.tree}, we have:

\begin{proposition} Let $\phi$ be a holomorphic function on $Tr_3$.
For any triangle $\Delta$, there is $\alpha_\Delta$, a bijection  $\{a, b\} \to \{M_1, M_2\}$, such that,
$\phi(c \Delta)= \alpha_\Delta(c)(\Delta)$, where $c $ is a letter $ \in \{a, b\}$. The same action is induced
on marked triangles, that is, $\phi(c \overline{\Delta}) = \alpha_\Delta(c) (\overline{\Delta})$.

\end{proposition}

 \subsubsection{Special holomorphic functions.} As in the case of $T_3$, some 
holomorphic functions appear special.
They are obtained, in a deteministic way, that is all the $\alpha_\Delta$ are the same, for example,
$\alpha_\Delta(a)= M_1, \alpha_\Delta(b) = M_2$ (for any $\Delta$). Note however that this depends on the choice of $M_1 $ and $M_2$,  or equivalently that of $\psi_1$ and $\psi_2$ (\S \ref{rough})

 \subsubsection{Action of $Aut(Tr_3)$.}  Let $Aut_0$ denote the stabilizer 
of $\Delta_0$ in $Aut(Tr_3)$,  the automorphism group of $Tr_3$. 
Applying elements of $Aut$ allows one to restore all
randomness  in the construction of holomorphic functions, that is the choice of the
$\alpha_\Delta$. From it one infers that, two holomorphic 
functions taking the same values
on $\Delta_0$ are equivalent up to composition in $Aut_0$.

Now, if one considers the full  action of $Aut(Tr_3)$, one gets:

\begin{fact} Two holomorphic functions $\phi$ and $\phi^\prime$ are equivalent
in $Aut(Tr_3)$, iff, they take the same values on two triangles $\Delta$
and $\Delta^\prime$

\end{fact}

 \subsubsection{Space of triangles. Orbital structure.} 
The space of marked (oriented) triangles
is  ${\mathcal T}=  \{(p, e, f) \in {\Bbb C}^3, (e, f) \neq (0, 0)\}$ 
(as was said previously, degenerate triangles, i.e. with $e=0$ or $f= 0$, are also considered).  We have seen, how holomorphic
functions induce, a kind of random dynamics on ${\mathcal T}$. Consider the relation $\sim$ defined by
$${\bf \Delta} \sim {\bf \Delta}^\prime \iff \exists \phi\; \mbox{holomorphic}, \Delta, \Delta^\prime\; \hbox{triangles of}\; Tr_3,  {\bf \Delta}= \phi (\Delta), {\bf \Delta}^\prime = \phi(\Delta^\prime)$$

In other words, thanks to the Fact above: ${\bf \Delta} \sim {\bf \Delta}^\prime $, iff, there is $\phi$
holomorphic, $\gamma \in Aut(Tr_3)$, such that ${\bf \Delta} $ and  $ {\bf \Delta}^\prime $ are images
of (the central triangle) $\Delta_0$ by $\phi$ and $\phi \circ \gamma$, respectively.

Previous developments show that $\sim$ is a well defined {\it equivalence} relation on ${\mathcal T}$
(despite randomness, and arbitrary choose of $\psi_1, \psi_2$...).

It is worthwhile to investigate ergodic theory of this relation!

\subsubsection{Action of the similarity group, Conformal triangles}
All involved dynamical notions  $I, M_1$, $M_2$, $\sim$...  commute with the action of the similarity group  $SG$ of ${\Bbb C}$. The marked vertex of the triangle can be identified under this action with
$0 \in {\Bbb C}$. It remains to consider pairs $(e, f) \in {\Bbb C}^2-(0,0)$,  up to homothety. This is exactly
${\Bbb C}P^1$, which is interpreted as the space of conformal triangles (allowed to be degenerate).

In particular, $\sim$ passes to  an equivalence relation $ \approx$ on ${\Bbb C}P^1$.

Let us observe here that an adapted parallel construction in the case of $T_3$, yields an equivalence
relation $\approx$ which is trivial. Here, our $\approx$ is far from being so:  not only the dynamics on triangles is not trivial, but even their conformal type is strongly transformed.

\subsection{Language of correspondences}

Denote  $(x, y, z)= (p-u, -u, -u+v)$. Then, $(p, u, v)= (x+y, -y,  -y+z)$.
Therefore, the ``graph'' of the multi- valued transformation $M$, is the subset of
$ (p, e, f; x, y, z) \in {\Bbb C}^3 \times {\Bbb C}^3$, such that:
\begin{eqnarray}
e+f-y+ ( -y+z) &=& 0 \nonumber \\
e^2+f^2+y^2+(-y+z)^2 &=& 0 \nonumber \\
x - (p+y) &=& 0 \nonumber
\end{eqnarray}


The space of solutions ${\mathcal S}$, is a 3-dimensional quadric contained in a 4-dimensional linear subspace. We have surjective induced projection $\pi_1: {\mathcal S} \to {\Bbb C}^3$ and
$\pi_2: {\mathcal S} \to {\Bbb C}^3$.  (Actually, instead of ${\Bbb C}^3$, one may better consider
the space of triangles ${\mathcal T}$, and ${\mathcal S} \subset {\mathcal T} \times {\mathcal T}$). The projection
$\pi_1$ and $\pi_2$ are branched converings. With all these objects,  we are bringing
out a natural situation  leading to
the useful notion of correspondences  (see \cite{Bul} for fundamental references). The point is to consider
the multi-valued dynamics: $ {\bf \Delta} \to \pi_2 (\pi_1^{-1} ({\bf \Delta}))$.


 \subsubsection{Projective version.} 
 As above, one can see what happens for conformal
 triangles by taking quotient under
the  action of the similarity group. More precisely, 
we consider the quotient of
$({\Bbb C}\times ({\Bbb C}^2-0)) \times ({\Bbb C}\times ({\Bbb C}^2-0))$
by the product action of $SG \times SG$ (recall that  $SG$ 
is the group of similarities
of ${\Bbb C}$).  We get ${\Bbb C}P^1 \times {\Bbb C}P^1$, 
and inside it, corresponding to
${\mathcal S}$,
a quadric ${\mathcal C}$, which thus determines 
a correspondence of ${\Bbb C}P^1$.

\subsubsection{Growth, Simulation}

As was said, we restrict ourselves in this case to formulation 
rather than a systematic
study of the dynamics. In order to have an idea, 
one can implement the multi-valuated
mapping :
$M(p, e, f)~=~$ \\
$ (\frac{ 2p +(e+f) -\sqrt{-3(e^2+f^2)-2ef} }{2},  
\frac{ (e+f)  - \sqrt{-3(e^2+f^2)-2ef} }{2},
 - \sqrt{-3(e^2+f^2)-2ef} )$
and gets, as example the following two pictures :

\centerline {\psfig{figure=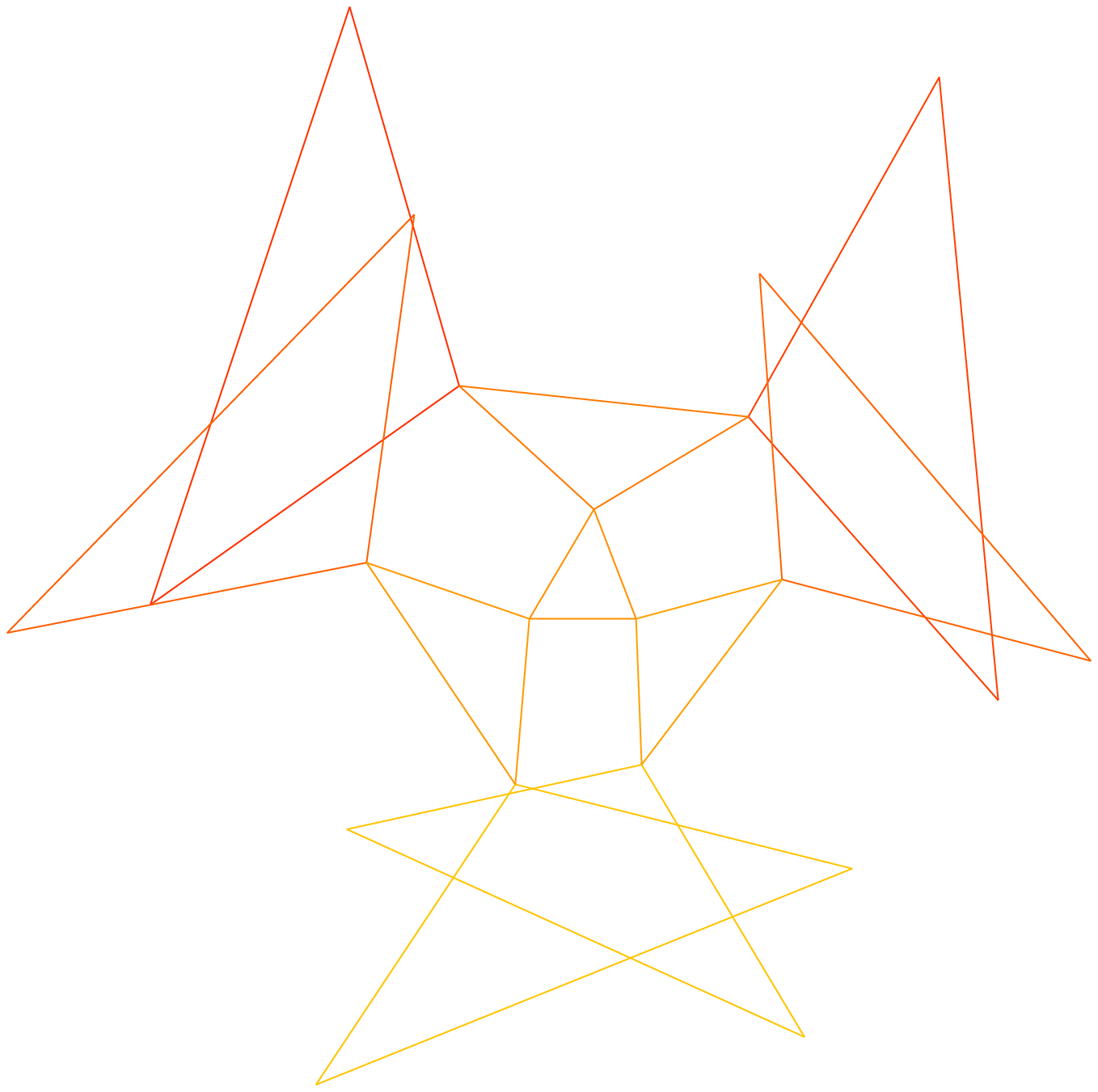,width=11cm,height=7cm} }
\centerline {\psfig{figure=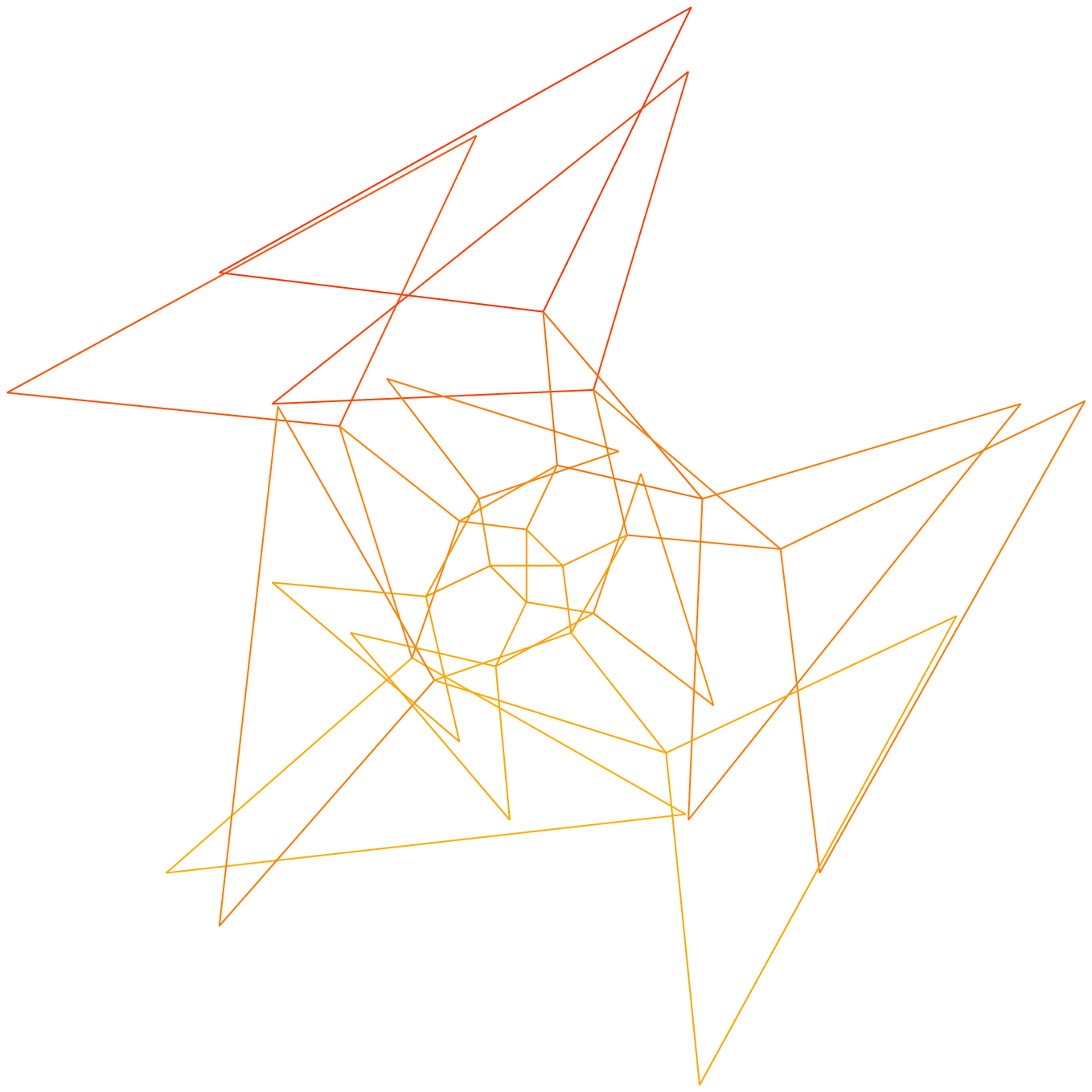,width=11cm,height=7cm}}
\begin{center}

{ Figure 8 : Images of balls in $Tr_3$ of radius  3 and 5 
by holomophic functions}

\end{center}






 \section{Remarks on other graphs} \label{remarks.other}

 \subsubsection{Graphs of valency 3.}  Any graph $X$ of valency 3 is (universally) 
covered by
$T_3$. Therefore any  holomorphic function
on $X$ gives rise to a holomorphic function on $T_3$, which is 
therefore fully understood
thanks to Theorem \ref{theorem.tr3}. We straightforwardly deduce from it the following:

\begin{theorem} Let $X$ be a graph of valency 3. If $X$ is the hexagonal tiling, 
then any holomorphic function $X \to
{\Bbb C}$ is the restriction of a similarity.

If $X$ has a cycle of length $\leq 5$, then, 
any holomorphic function on $X$ is constant
(here we tacitly assume $X$ is connected). In fact, 
in general, holomorphic functions are
trivial unless $X$ covers the hexagonal tiling.

\end{theorem}

As example, the cayley graph of $PSL(2, {\Bbb Z})$ has no non-constant  holomorphic functions.

$$\epsh{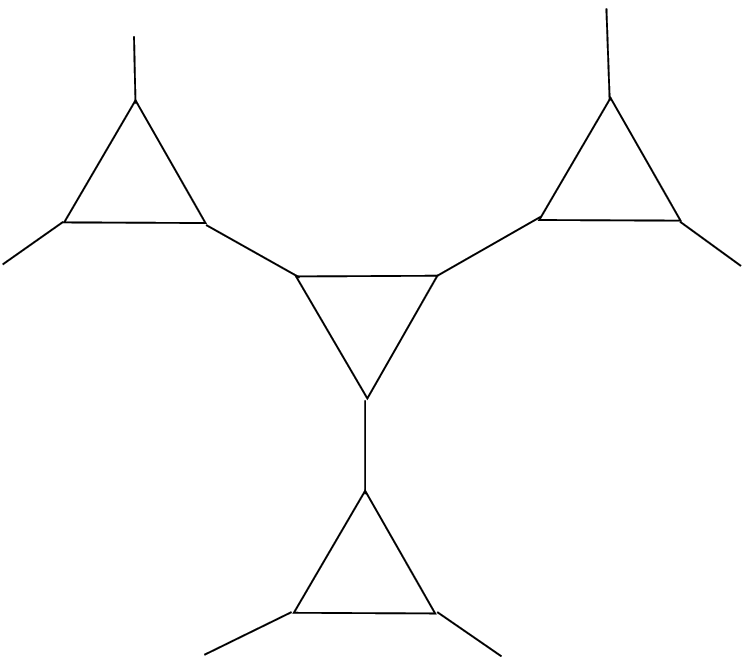}{2.5cm}$$ 
\begin{center}

{ Figure 9 : Cayley graph of $PSL(2, {\Bbb Z})$}

\end{center}

 \subsubsection{Cayley graph of ${\Bbb Z}^n$.}
Consider the cayley graph $X_n$ of ${\Bbb Z}^n$ (associated to the canonical basis together with its opposite), that is, the standard lattice in ${\Bbb R}^n$.

For $n=1$, harmonic functions are affine  : $n\mapsto a+bn$ and
  holomorphic functions are  constant.

Identify  $X_2$ with $\Bbb Z^2=\Bbb Z+i\Bbb Z$
embeded in $\Bbb C$.
One  easily   verifies that,  only if $k\leq 3$,
 the function $z\mapsto z^k$ is
harmonic on  $X_2$. On the other hand , the mappings
 $z\mapsto az+b$, and similarly $z \mapsto a\bar{z}+b$, are    holomorphic.  It is
very likely that these are all holomorphic functions on $X_2$?




\subsection{Non-rigidity, Trees of higher valency }


For a  tree $T$  of valency $\geq 4$, and two adjacent
vertices $s_0$ and $s_1$, on which a function $\phi$ takes values $z_0$
and $z_1$, there is at least 1 degree   of freedom
in extending holomorphically, progressively $\phi$ to neighbouring vertices. Assume, to simplify notation,  the valency equals $4$ and let $s_2, s_3$ and $s_4$ the other neighbours
of $s_0$.  The associated values $z_i = \phi(s_i)$, for $i= 2, 3, 4$, satisfy:
\begin{eqnarray}
\sum_{i=2}^{i=4} \delta_i    &=&  - \delta_1 \nonumber \\
 \sum_{i=2}^{i=4} \delta_i^2 &=& - \delta_1^2 \nonumber
 \end{eqnarray}

where $\delta_i = z_i-z_0$. Only $\delta_1$ is given, thus,  we have 2 equations and 3 unknowns, and hence,  a 1-dimensional complex space of solutions.

 \subsubsection{Bounded holomorphic functions.} We will now observe that with this freedom, we can choose $\phi$ to be bounded (and non-constant). (Even not formally stated previously,  this is not possible for $T_3$ and $Tr_3$).    To do this, we show that, there is
a positive real number  $r <1$, such that, for any $\delta_1$,  the above
equations admit solutions satisfying $ \vert \delta_i \vert \leq r \vert  \delta_1 \vert$,
for $i = 2, 3, 4$. (With this, one constructs a holomorphic function $\phi$
having an exponentially  decreasing oscillation along edges, and hence  bounded).
By homogeneity, it suffices to consider the case $\delta_1= -1$. To get a solution in this
case,  we write $\delta_2= r_1e^{i \theta}$, $\delta_3 = \bar{\delta_2} = r_1e^{-i \theta}$
and $\delta_4 = r_2$ real and positive. We claim that equations:
$$2r_1 \cos \theta +r_2 = 1, \; 2r_1^2 \cos 2 \theta + r_2^2 = -1$$
admit  solutions with $0 < r_1, r_2 <1$.  This is done by elementary analysis around
$r_1 = r_2 = 1$, and $\theta =\frac{ \pi}{2}$.

\begin{remark}
{\em As was suggested in \S \ref{graph.case}, a theory parallel to that
of $T_3$, may be developed for a tree $T_N$ of valency $N$, if one considers
``$(N-1)$-holomorphic functions''. This would surely give  beautiful pictures.
 }

 \end{remark}

\section{The conjugate part Problem} \label{conjugate.problem}

Let us recall the problem.  We consider a {\it real}  
function $f$ on a graph $X$, and we ask, when (and how) it has
 a conjugate part, that is another real function $g$, 
such that $\phi= f+ig$ is holomorphic?
 Actually, we have finally give up to systematically  
investigate this problem in the present article, 
since of unexpected subtleties. Only partial results will be exposed here.
To mention a rough conclusion,
  we say that also here a random dynamics is in order, and
as
 it is natural to expect,
the problem has some  cohomological flavour.

 \subsubsection{Conditions.} It turns out that the same conditions, stated in the Riemannian
case, as explained after Definition \ref{definition.holomorphic},   
hold in the discrete case.  For any vertex $s$,
$g$ satisfies:
 $$\parallel \nabla_s f \parallel =
\parallel \nabla_s g \parallel, \;
 \langle \nabla_s f, \nabla_s g \rangle = 0,\; \langle \nabla_s g, 
(1, \ldots, 1) \rangle= 0  $$

 The last condition expresses the fact that $g$ is harmonic.  
If $s$ has valency $n$,
 these conditions mean exactly:  $\nabla_s g $ must belong to $ {\mathcal C}^f_s$, 
a $(n-3)$-sphere of radius 
$ \parallel \nabla_s f \parallel $ in ${\Bbb R}^n$.

 From this ``freedom'', say when ${\mathcal C}^f$ has higher dimension, one may ask
 if $f$ is allowed to have more than one conjugate part (up to constant). 
It turns out this is
 possible, for instance in the case of the (homogeneous) tree of valency 4.  
Notice however, that
 because of the norm gradient equality, the space of conjugates (of a given $f$)
 is compact (up to constants).

 \subsubsection{Non-existence examples.} On the other hand, 
on the same tree (of valency 4), there are harmonic
  functions with no conjugate part.  For an example (see Figure 10\label{}), 
consider $A$ a vertex, $B, C, D$ and $E$
  its  neighbours. Consider the sub-tree ${\mathcal A}$ consisting of $A, B, C, D, E$ and
the  3 remaining  vertices  of each $B$ and $C$.

  Take $f$ to be everywhere 0 on ${\mathcal A}$,  except on    $D$ and $E$
   where it takes two opposite non vanishing
  values, say 1 and $-1$.

   Suppose a conjugate $g$ exists. Adding a constant, we can
   assume $g(A) = 0$. The norm gradient equality implies
   that,  like $f$, $g$ vanishes on all neighbours of $B$ and $C$.

   On the other hand, writing the conjugate part equations on $A$, we get
   solutions which must vanish on $D$ and $E$, but not on $B$ and $C$,
   which is a contradiction.
 Observe  here that $f$ can be easily extended (not uniquely) to a harmonic function
   on all the tree.

   One may believe  the ``degeneracy'' of $f$, that is the vanishing 
of many of its oscillations
   and gradients,  was the responsible on non-existence of its conjugate part. 
This is not true, all harmonic functions
   sufficiently near $f$ have no conjugate part as well. Indeed, having a 
conjugate part is a closed condition.  This follows from the 
compactness mentioned above (due  the
   the norm gradient equality).


$$\epsh{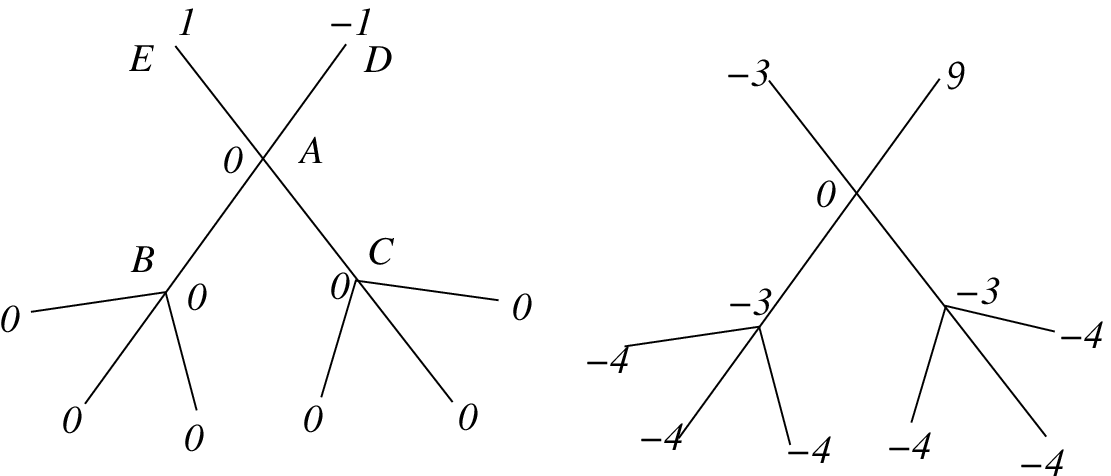}{3.5cm}$$ 

 \begin{center}
 { Figure 10 : No conjugate part}

 \end{center}

 \subsubsection{An existence result.}

  From the classification of holomorphic functions
  on $T_3$, one infers that  a
harmonic function has  conjugate part if
and only if  its  gradient has a   
constant norm (throughout $T_3$) (Actually, this corollary
    is equivalent to the classification itself). 
Our principal existence result is that
  this condition, norm constancy of the gradient, 
is sufficient for trees of higher
  valences. This is, by no means, a necessary condition, 
but seems nevertheless, natural,
  and gives an interesting  class of examples.


   \begin{theorem} Let $f$ be a harmonic on a (homogeneous) tree $X$,
  of (finite) valency $ \geq 4$, and assume its gradient has a constant norm,
  say, $\parallel  \nabla_s f  \parallel = 1$, $\forall s \in X$.
  Then $f$ has a conjugate part.  There is  $(\nu-3)$
   degrees of freedom (at each vertex)
   in choosing this conjugate part, where $\nu$ is the valency.   The space of
   conjugates of $f$ is compact.

   \end{theorem}


   \begin{proof} Choose $O \in X$, and for $r \in {\Bbb  N}$, 
let $B_r$ be the (closed) ball centred at $O$ of radius $r$
   (i.e. the set of vertices that can be reached from $O$ by a sequence of  at most
   $r$ edges).  Consider also the sphere $S_r$, the boundary of $B_r$.

   Trees are characterized by the following ``strict convexity'' property of $B_r$: 
if $s \in S_r$, then,
   only one neighbour $p(s) $ of $s$ lies in $B_r$, in fact in $B_r-S_r$.  
All the others are outside
   $B_r$, in fact they belong to $S_{r+1}$. Furthermore a vertex of $S_{r+1}$ 
is adjacent
   to exactly one vertex of  $S_r$.


   Assume $g$ is defined on $B_r$, and let us show how to extend (with some freedom)
   it to $B_{r+1}$.  Let $s \in S_r$. Then, $g$
   is already defined on exactly one neighbour $p(s)\in S_{r-1}$ of $s$.

   In order, to extend $g$ to other
   neighbours  ($\in S_{r+1}$), we just have to give $\nabla_s g$.   
Actually,  exactly one component, say $a_1$,
of   $\nabla_sg$ is already given, that corresponding to the edge $p(s)-s$.   
The gradient
     $\nabla_s g$  must belong to the sphere ${\mathcal C}^f_s$. 
The question reduces then to the fact that this sphere ${\mathcal C}^f_s$ 
contains a vector with first projection $a_1$ (again  to be precise,
     by first projection, we mean that corresponding to the edge $p(s)-s$). 
The proof of the theorem will be completed thanks to the following lemma. 
Its significance is that,
the spheres     ${\mathcal C}^f_s$ and ${\mathcal C}^f_{p(s)}$ have exactly the same
projection on the coordinate corresponding to the edge $p(s)-s$.
      \end{proof}

\begin{lemma}
Let $n \geq 4$, $\delta=(\delta_1,\ldots,\delta_n)\in \RR^n\setminus \lbrace 0\rbrace,
 \sum \delta_i=0$, $\sum \delta_i^2=1$.
     Consider the sphere
$${\mathcal C}=\lbrace a=(a_1,...,a_n)\in\RR^n,
\sum a_i \delta_i=0,\sum a_i^2=1,  \sum a_i=0 \rbrace$$
 Then
the image of  ${\mathcal C}$ by the  projection on the first
 factor of $\RR^n$ is the  segment
$[-\alpha,\alpha]$ where $\alpha=\sqrt{\frac{n-1}{n}-\delta_1^2}$.
      In particular,
     {\em this depends only on  $\delta_1$}.


\end{lemma}

  Before proving it let us make some remarks about this lemma. Firstly, it
    was stated in the case $\parallel \delta \parallel = 1$, but, obviously, 
a general formula is available
  for any constant norm.

  We also notice that a  similar result is true for $n=3$, the projection
in this case consists of two points : $\pm \sqrt{\frac{2}{3}-\delta_1^2}$. 
The theorem
  itself is true for $T_3$, in fact as was said above, it is optimal in this case. 
Observe however, that the theorem does not extend to the (trivial) case of 
the 2-valenced tree
  $T_2 (\approx {\Bbb Z})$.

  Finally,  observe that,
our second example   of a  harmonic function  without conjugate part 
(see Fig. 10)
corresponds
to  the case where $\alpha=0$  (i.e. $\delta_1=\sqrt 3 /2$ and 
$\delta_2=\delta_3=\delta_4=-\sqrt 3 /6$).

   \begin{proof}
Consider the mapping
   $\pi : a\mapsto a_1$. The projection of the gradient of $\pi$
on the hyperplane
$H=\lbrace (x_1,...,x_n)\in \RR^n,~ \sum x_i=0\rbrace$ is collinear
to $v=(1-n,1,...,1)$, and so if  $a$ is a critical point of  $\pi$ restricted to
the submanifold  ${\mathcal C}$, then there exist $\alpha, \beta$
so that  $a=\alpha v+\beta \delta$.
>From  the fact that $a$ is
 orthogonal  to  $\delta$,  we infer that : $\beta=-\alpha <v,\delta>$.
But, since $ \sum \delta_i=0$, we have $ <v,\delta>= -n\delta_1$.
Finally, we have $a=\alpha( (1-n)+n \delta_1^2, ~1+n\delta_1\delta_2,...,
1+n\delta_1\delta_n)$. Thus,
 $\parallel v-<v,\delta> \delta \parallel^2=n(n-1-n\delta_1^2)$, and hence,
 $a_1^2=\frac{n-1}{n}-\delta_1^2$.
 \end{proof}




\begin{remark}
{\em
(Canonical conjugate?)  The conjugate part problem is sub-determined: solutions
   are not unique (when they exist). It is interesting to  define  ways of selecting unique solutions.
   Maybe, by means of a companion variational problem?
   }
  \end{remark}




${}$

\end{document}